# A PARTICULAR COLLECTION OF SERVICE TIME DISTRIBUTIONS PARAMETERS STUDY AND IMPACT IN SOME M|G|∞ SYSTEM BUSY PERIOD AND BUSY CYCLE PARAMETERS

Manuel Alberto M. Ferreira

Instituto Universitário de Lisboa (ISCTE-IUL), Business Research Unit (BRU-IUL) and Information Sciences, Technologies and Architecture Research Center (ISTAR-IUL), **PORTUGAL**

E-mail : manuel.ferreira@iscte.pt

**ABSTRACT**

*The problems arising when the moments of service time distributions, for which the M|G|∞ queue system busy period and busy cycle become very easy to study, are presented and it is shown how to overcome them. The busy cycle renewal function and the "peakedness" and the "modified peakedness" for the M|G|∞ busy period and busy cycle in the case of those service time distributions are also computed.*

***Keywords*:** *Service time, collection, distributions, moments, M|G|∞ queue.*

**INTRODUCTION**

When, in the M|G|∞ queue system, the service time length is a random variable with a distribution function belonging to the collection

$$G(t) = 1 - \frac{\left(1 - e^{-\rho}\right)\left(\lambda + \frac{\lambda p + \beta}{1-p}\right)}{\lambda e^{-\rho}\left(e^{\left(\lambda + \frac{\lambda p + \beta}{1-p}\right)t} - 1\right) + \lambda}, t \geq 0, -\lambda \leq \beta \leq \frac{\lambda\left(1 - pe^{\rho}\right)}{e^{\rho} - 1}, 0 \leq p < 1 \quad (1.1),$$

the busy period length probability distribution is exponential with an atom at the origin and the busy cycle length probability distribution is the mixture of two exponential distributions, see Ferreira (2005) and (Ferreira and Andrade, 2009). But although it is so easy to study the busy period and the busy cycle in this situation it is very difficult to compute the service time moments.
    Some results, precisely about the moment's computation of random variables with distribution functions given by this collection are given.
    In the end are presented formulae that give the busy cycle renewal function and the "peakedness" and the "modified peakedness" to the busy period and the busy cycle of the M|G|∞ system for those service time distributions, see Ferreira (2004, 2013, 2013a).

This work is built on the presented in Ferreira (2007) which is so corrected, generalized and updated.

**MOMENTS COMPUTATION**

Be $G(t), t \geq 0$ a distribution function and $g(t) = \dfrac{dG(t)}{dt}$.

The differential equation $(1-p)\dfrac{g(t)}{1-G(t)} - \lambda p - \lambda(1-p)G(t) = \beta$, where $\lambda > 0$ and $-\lambda \leq \beta \leq \dfrac{\lambda(1-pe^{\rho})}{e^{\rho}-1}$, $0 \leq p < 1$ ($\rho = \lambda\alpha$, being $\alpha$ the mean of $G(t)$) has (1.1) as solution (see Ferreira (2005)).

If, in (1.1), $G_i(t)$ is the solution associated to $\rho_i$, $i = 1,2,3,4$ it is easy to see that

$$\frac{G_4(t) - G_2(t)}{G_4(t) - G_1(t)} \cdot \frac{G_3(t) - G_1(t)}{G_3(t) - G_2(t)} = \frac{e^{-\rho_4} - e^{-\rho_2}}{e^{-\rho_4} - e^{-\rho_1}} \cdot \frac{e^{-\rho_3} - e^{-\rho_1}}{e^{-\rho_3} - e^{-\rho_2}} \qquad (2.1)$$

as it had to happen since it is a Riccati equation.

And computing,

$$\int_0^\infty [1 - G(t)]dt = \int_0^\infty \frac{(1-e^{-\rho})\left(\lambda + \dfrac{\lambda p + \beta}{1-p}\right)}{\lambda e^{-\rho}\left(e^{\left(\lambda + \frac{\lambda p + \beta}{1-p}\right)t} - 1\right) + \lambda} dt =$$

$$= \frac{(1-e^{-\rho})\left(\lambda + \dfrac{\lambda p + \beta}{1-p}\right)}{\lambda} \int_0^\infty \frac{1}{e^{-\rho}\left(e^{\left(\lambda + \frac{\lambda p + \beta}{1-p}\right)t} - 1\right) + 1} dt =$$

$$= \frac{(1-e^{-\rho})\left(\lambda + \dfrac{\lambda p + \beta}{1-p}\right)}{\lambda} \int_0^\infty \frac{e^{-\left(\lambda + \frac{\lambda p + \beta}{1-p}\right)t}}{e^{-\rho} - e^{-\rho}e^{-\left(\lambda + \frac{\lambda p + \beta}{1-p}\right)t} + e^{-\left(\lambda + \frac{\lambda p + \beta}{1-p}\right)t}} dt =$$

$$= \frac{(1-e^{-\rho})\left(\lambda + \dfrac{\lambda p + \beta}{1-p}\right)}{\lambda} \int_0^\infty \frac{e^{-\left(\lambda + \frac{\lambda p + \beta}{1-p}\right)t}}{e^{-\rho} + (1 - e^{-\rho})e^{-\left(\lambda + \frac{\lambda p + \beta}{1-p}\right)t}} dt =$$

$$= \frac{(1-e^{-\rho})\left(\lambda + \frac{\lambda p + \beta}{1-p}\right)}{\lambda} \cdot \frac{-1}{(1-e^{-\rho})\left(\lambda + \frac{\lambda p + \beta}{1-p}\right)} \cdot \left[\log\left(e^{-\rho} + (1-e^{-\rho})e^{-\left(\lambda + \frac{\lambda p + \beta}{1-p}\right)t}\right)\right]_0^\infty =$$

$$= -\frac{1}{\lambda}(\log e^{-\rho} - \log 1) = \frac{-\rho}{-\lambda} = \alpha$$

as it had to be because are considered positive random variable.

The density associated to $G(t)$ given by (1.1) is

$$g(t) = \frac{(1-e^{-\rho})e^{-\rho}\left(\lambda + \frac{\lambda p + \beta}{1-p}\right)^2 e^{-\left(\lambda + \frac{\lambda p + \beta}{1-p}\right)t}}{\lambda \left[e^{-\rho} + (1-e^{-\rho})e^{-\left(\lambda + \frac{\lambda p + \beta}{1-p}\right)t}\right]^2}, t > 0, -\lambda \leq \beta \leq \frac{\lambda(1-pe^{-\rho})}{e^\rho - 1}, 0 \leq p < 1 \quad (2.2).$$

So,

$$\int_0^\infty t^n g(t) dt = \frac{(1-e^{-\rho})e^{-\rho}\left(\lambda + \frac{\lambda p + \beta}{1-p}\right)^2}{\lambda} \cdot \int_0^\infty t^n \frac{e^{-\left(\lambda + \frac{\lambda p + \beta}{1-p}\right)t}}{\left[e^{-\rho} + (1-e^{-\rho})e^{-\left(\lambda + \frac{\lambda p + \beta}{1-p}\right)t}\right]^2} dt \cdot$$

But, $\int_0^\infty t^n \frac{e^{-\left(\lambda + \frac{\lambda p + \beta}{1-p}\right)t}}{\left[e^{-\rho} + (1-e^{-\rho})e^{-\left(\lambda + \frac{\lambda p + \beta}{1-p}\right)t}\right]^2} dt \geq \int_0^\infty t^n e^{-\left(\lambda + \frac{\lambda p + \beta}{1-p}\right)t} dt =$

$$= \frac{1}{\lambda + \frac{\lambda p + \beta}{1-p}} \frac{n!}{\left(\lambda + \frac{\lambda p + \beta}{1-p}\right)^n}, \beta \neq -\lambda.$$

And, $\int_0^\infty t^n \frac{e^{-\left(\lambda + \frac{\lambda p + \beta}{1-p}\right)t}}{\left[e^{-\rho} + (1-e^{-\rho})e^{-\left(\lambda + \frac{\lambda p + \beta}{1-p}\right)t}\right]^2} dt \leq e^{2\rho} \int_0^\infty t^n e^{-\left(\lambda + \frac{\lambda p + \beta}{1-p}\right)t} dt =$

$$= \frac{e^{2\rho}}{\lambda + \frac{\lambda p + \beta}{1-p}} \frac{n!}{\left(\lambda + \frac{\lambda p + \beta}{1-p}\right)^n}, \beta \neq -\lambda.$$

So, calling T the random variable corresponding to $G(t)$:

$$\frac{(1-e^{-\rho})e^{-\rho}}{\lambda} \frac{n!}{\left(\lambda + \frac{\lambda p + \beta}{1-p}\right)^{n-1}} \leq E[T^n] \leq \frac{e^\rho - 1}{\lambda} \frac{n!}{\left(\lambda + \frac{\lambda p + \beta}{1-p}\right)^{n-1}},$$

$$, -\lambda < \beta \leq \frac{\lambda(1-pe^{-\rho})}{e^\rho - 1}, 0 \leq p < 1, n = 1, 2, \ldots \quad (2.3).$$

**Notes:**
- The expression (2.3), giving bounds for $E[T^n]$, guarantees its existence,
- For $n=1$ the expression (2.3) is useless since $E[T] = \alpha$. Note, curiously, that the upper bound is $\frac{e^\rho - 1}{\lambda}$, the M|G|∞ system busy period mean value,
- For $n=2$, subtracting to both bounds $\alpha^2$, it is possibly get from expression (2.3) bounds for $VAR[T]$,
- For $\beta = -\lambda, E[T^n] = 0, n = 1, 2, \ldots$, evidently.

See, however, that (1.1) can be written like:

$$G(t) = \frac{1 + \frac{\lambda p + \beta}{1-p}(1-e^\rho)e^{-\left(\lambda + \frac{\lambda p + \beta}{1-p}\right)t}}{1 - (1-e^\rho)e^{-\left(\lambda + \frac{\lambda p + \beta}{1-p}\right)t}}, t \geq 0, -\lambda \leq \beta \leq \frac{\lambda(1-pe^\rho)}{e^\rho - 1}, 0 \leq p < 1 \quad (2.4)$$

and, for $\rho < \log 2$,

$$G(t) = \left(1 + \frac{\frac{\lambda p + \beta}{1-p}}{\lambda}(1-e^\rho)e^{-\left(\lambda + \frac{\lambda p + \beta}{1-p}\right)t}\right) \cdot \sum_{k=0}^{\infty} (1-e^\rho)^k e^{-k\left(\lambda + \frac{\lambda p + \beta}{1-p}\right)t},$$

$$, t \geq 0, -\lambda \leq \beta \leq \frac{\lambda(1-pe^\rho)}{e^\rho - 1}, 0 \leq p < 1 \quad (2.5).$$

After (2.5) it is easy to derive the T Laplace Transform for $\rho < \log 2$. And, so,

- For $\rho < \log 2$

$$E[T^n] = -\left(1 + \frac{\frac{\lambda p + \beta}{1-p}}{\lambda}\right) n! \sum_{k=1}^{\infty} \frac{(1-e^\rho)^k}{k\left(\lambda + \frac{\lambda p + \beta}{1-p}\right)^n}, -\lambda < \beta \leq \frac{\lambda(1-pe^\rho)}{e^\rho - 1}, 0 \leq p < 1,$$

$$, n = 1, 2, \ldots \qquad (2.6).$$

**Notes:**

- $E[T] = -\left(1 + \frac{\frac{\lambda p + \beta}{1-p}}{\lambda}\right) \sum_{k=1}^{\infty} \frac{(1-e^\rho)^k}{k\left(\lambda + \frac{\lambda p + \beta}{1-p}\right)} = \frac{1}{\lambda} \sum_{k=1}^{\infty} (-1)^{k+1} \frac{(1-e^\rho)}{k} =$

$$= \frac{1}{\lambda} \log e^\rho = \frac{\rho}{\lambda} = \alpha.$$

- For $n \geq 2$ only a finite number of parcels can be considered in the infinite sum. Calling M this number, to get an error lesser than $\varepsilon$ it must be fulfilled simultaneously

a) $M > \dfrac{1}{\lambda + \dfrac{\lambda p + \beta}{1-p}} - 1,$

b) $M > \log_{(e^\rho - 1)} \dfrac{\varepsilon e^\rho \lambda}{n!\left(\lambda + \dfrac{\lambda p + \beta}{1-p}\right)} - 1.$

So, it is evident now that this distributions collection moment's computation is a complex task. This was already true for the study of Ferreira (1998) where the results presented are a particular situation of these ones for $p = 0$.

To consider the approximation

$$E_m^n = \sum_{k=1}^{\infty} \left(\frac{k}{m}\right)^n \left[G\left(\frac{k}{m}\right) - G\left(\frac{k-1}{m}\right)\right], -\lambda < \beta \leq \frac{\lambda(1-pe^\rho)}{e^\rho - 1}, 0 \leq p < 1, n = 1, 2, \ldots \qquad (2.7)$$

may be helpful since $\lim_{m \to \infty} E_m^n = E[T^n], n = 1, 2, \ldots$ (Ferreira and Andrade, 2012c) that allow the moments numerical computation.

### BUSY CYCLE RENEWAL FUNCTION COMPUTATION

The busy cycle (an idle period followed by a busy period) renewal function value of the M|G|∞ queue, at $t$, gives the mean number of busy periods that begin in $[0, t]$, see

Ferreira (2004). If the service time is a random variable with distribution function given by a member of the collection (1.1), calling the value of the renewal function at $t$ $R(t)$:

$$R(t) = e^{-\rho}(1+\lambda t) + (1-e^{-\rho})\frac{\lambda p + \beta}{\lambda + \beta} e^{-\left(\lambda + \frac{\lambda p + \beta}{1-p}\right)t} + (1-e^{-\rho})\frac{\lambda p + \beta}{\lambda + \beta}, -\lambda \leq$$

$$\leq \beta \leq \frac{\lambda(1-pe^{\rho})}{e^{\rho}-1}, 0 \leq p < 1 \qquad (3.1).$$

For $p = 0$ it is obtained the result presented in Ferreira (2004).

### THE "PEAKEDNESS" AND THE "MODIFIED PEAKEDNESS" FOR THE M|G|∞ QUEUE BUSY PERIOD AND BUSY CYCLE

The M|G|∞ queue busy period "peakedness" is the Laplace Transform of its length at $1/\alpha$, Ferreira (2013, 2013a). It is a parameter that characterizes the busy period distribution length and contains information about all its moments. For the collection of service distributions (1.1) the "peakedness", named $pi$, is

$$pi = \frac{e^{-\rho}(\lambda + \beta)(\rho + 1) - \lambda p - \beta}{\lambda(e^{-\rho}(\rho + \alpha\beta) + 1 - p)}, -\lambda \leq \beta \leq \frac{\lambda(1-pe^{\rho})}{e^{\rho}-1}, 0 \leq p < 1 \qquad (4.1).$$

In Ferreira (2013, 2013a) is also introduced another measure, the "modified peakedness" got after the "peak" taking out the terms that are permanent for the busy period in different service distributions and putting over the common part. Calling it $qi$:

$$qi = pi\frac{\rho}{e^{\rho}-\rho-1} + 1$$

and so, for the distributions given by collection (1.1):

$$qi = \frac{e^{-\rho}(\lambda + \beta)(\rho + 1) - \lambda p + \beta}{\lambda(e^{-\rho}(\rho + \alpha\beta) + 1 - p)} \frac{\rho}{e^{\rho}-\rho-1} + 1, -\lambda \leq \beta \leq \frac{\lambda(1-pe^{\rho})}{e^{\rho}-1}, \ 0 \leq p < 1 \quad (4.2).$$

For the busy cycle of the M|G|∞ queue, analogously, it may be defined the "peakedness", Ferreira (2013a)), now called $pi'$, and for the service distributions given by the collection (1.1) it is

$$pi' = \alpha\frac{e^{-\rho}(\lambda + \beta)(\rho + 1) - \lambda p - \beta}{(\rho + 1)(e^{-\rho}(\rho + \alpha\beta) + 1 - p)}, -\lambda \leq \beta \leq \frac{\lambda(1-pe^{\rho})}{e^{\rho}-1}, 0 \leq p < 1 \qquad (4.3)$$

and the "modified peakedness", that now called $qi'$, given by $pi'\frac{\rho}{e^{\rho}-\rho} + 1$, and for the service distributions given by the collection (1.1) it is

$$qi' = \alpha \frac{e^{-\rho}(\lambda+\beta)(\rho+1)-\lambda p-\beta}{(\rho+1)(e^{-\rho}(\rho+\alpha\beta)+1-p)}\frac{\rho}{e^{\rho}-\rho}+1, -\lambda \leq \beta \leq \frac{\lambda(1-pe^{\rho})}{e^{\rho}-1}, 0 \leq p < 1 \quad (4.4).$$

**REFERENCES**

1) ANDRADE, M. (2010). A note on foundations of probability. *Journal of Mathematics and Technology*, 1(1), 96-98.

2) CARRILLO, M. J. (1991). Extensions of Palm's theorem: a review. *Management Science*, 37(6), 739-744.

3) FERREIRA, M. A. M., RAMALHOTO, M. F. (1994). Estudo dos parâmetros básicos do período de ocupação da fila de espera $M|G|\infty$. *A Estatística e o Futuro e o Futuro da Estatística*. Actas do I Congresso Anual da S.P.E., Edições Salamandra, Lisboa, 47-60.

4) FERREIRA, M. A. M. (1998). Momentos de Variáveis Aleatórias com Função de Distribuição dadas pela Colecção $G(t)=1-\frac{(1-e^{-\rho})(\lambda+\beta)}{\lambda e^{-\rho}(e^{(\lambda+\beta)}-1)+\lambda}, t \geq 0, \quad -\lambda \leq \beta \leq \frac{\lambda}{e^{\rho}-1}$. Revista Portuguesa de Gestão, II. ISCTE, Lisboa, 67-69.

5) FERREIRA, M. A. M. (1998a). Application of Riccati equation to the busy period study of the $M|G|\infty$ system. *Statistical Review*, 1st Quadrimester, INE, 23-28.

6) FERREIRA, M. A. M. (1998b). Computational simulation of infinite servers systems. *Statistical Review*, 3rd Quadrimester, INE, 23-28.

7) FERREIRA, M. A. M. (2001). $M|G|\infty$ queue heavy-traffic situation busy period length distribution (power and Pareto service distributions). *Statistical Review*, 1st Quadrimester, INE, 27-36.

8) FERREIRA, M. A. M. (2002). The exponentiality of the $M|G|\infty$ queue busy period. *Actas das XII Jornadas Luso-Espanholas de Gestão Científica*, Volume VIII-Economia da Empresa e Matemática Aplicada. UBI, Covilhã, Portugal, 267-272.

9) FERREIRA, M. A. M. (2004), *M|G|∞ Queue Busy Cycle Renewal Function for Some Particular Service Time Distributions. Proceedings of Quantitative Methods in Economics (Multiple Criteria Decision Making XII)*. Virt, Slovakia, 42-47.

10) FERREIRA, M. A. M. (2005). Differential equations important in the $M|G|\infty$ queue system transient behavior and busy period study.

*Proceedings of 4th International Conference APLIMAT 2005*, Bratislava, Slovakia, 119-132.

11) FERREIRA, M. A. M. (2005a), A Equação de Riccati no Estudo da Ocupação de um Sistema M|G|∞ (Uma Generalização). *Actas do I Congresso de Estatística e Investigação Operacional da Galiza e Norte de Portugal/VII Congreso Galego de Estatística e Investigación de Operacións*. Guimarães, Portugal.

12) FERREIRA, M. A. M. (2007). M|G|∞ system parameters for a particular collection of service time distributions. *Proceedings of 6th International Conference APLIMAT 2007*, Bratislava, Slovakia, 131-137.

13) FERREIRA, M. A. M., ANDRADE, M., FILIPE, J. A. (2009). Networks of queues with infinite servers in each node applied to the management of a two echelons repair system. *China-USA Business Review*, 8(8), 39-45 and 62.

14) FERREIRA, M. A. M., ANDRADE, M. (2009). The ties between the $M|G|\infty$ queue system transient behavior and the busy period. *International Journal of Academic Research*, 1(1), 84-92.

15) FERREIRA, M. A. M., ANDRADE, M. (2010). Looking to a $M|G|\infty$ system occupation through a Riccati equation. *Journal of Mathematics and Technology*, 1 (2), 58-62.

16) FERREIRA, M. A. M., ANDRADE, M. (2010a). $M|G|\infty$ system transient behavior with time origin at the beginning of a busy period mean and variance. *APLIMAT-Journal of Applied Mathematics*, 3(3), 213-221.

17) FERREIRA, M. A. M., ANDRADE, M. (2011).Fundaments of theory of queues. *International Journal of Academic Research*, 3(1), part II, 427-429.

18) FERREIRA, M. A. M., ANDRADE, M. (2012). Busy period and busy cycle distributions and parameters for a particular $M|G|\infty$ queue system. *American Journal of Mathematics and Statistics*, 2(2), 10-15.

http://article.sapub.org/10.5923.j.ajms.20120202.03.html

19) FERREIRA, M. A. M., ANDRADE, M. (2012a). Queue networks models with more general arrival rates. *International Journal of Academic Research*, 4(1), part A, 5-11.

20) FERREIRA, M. A. M., ANDRADE, M. (2012b). Transient behavior of the $M|G|$m and $M|G|\infty$ system. *International Journal of Academic Research*, 4(3), part A, 24-33.

21) FERREIRA, M. A. M., ANDRADE, M. (2012c). A methodological note on the study of queuing networks. *Journal of Mathematics and Technology*, 3(1), 4-6.

22) FERREIRA, M. A. M., ANDRADE, M., FILIPE, J. A (2012). The age or excess of the M/G/∞ queue busy cycle mean value. *Computer and Information Science*, 5(5), 93-97. http://dx.doi.org/10.5539/cis.v5n5p93

23) FERREIRA, M. A. M. (2013). A $M|G|\infty$ queue busy period distribution characterizing parameter. *Computer and Information Science*, 6 (1), 83-88. http://dx.doi.org/10.5539/cis.v6n1p83

24) FERREIRA, M. A. M. (2013a). The modified peakedness as a M/G/∞ busy cycle distribution characterizing parameter. *International Journal of Academic Research*, 5(2), part A, 5-8.

25) FERREIRA, M. A. M. (2016). Results and applications in statistical queuing theory. *APLIMAT 2016 - 15th Conference on Applied Mathematics. Proceedings*. Bratislava, Slovakia, 362-375.

26) FIGUEIRA, J., FERREIRA. M. A. M. (1999). Representation of a pensions fund by a stochastic network with two nodes: an exercise. *Portuguese Revue of Financial Markets*, 1(3).

27) HERSHEY, J. C., WEISS, E. N., MORRIS, A. C. (1981). A stochastic service network model with application to hospital facilities. *Operations Research*, 29(1), 1-22, 1981.

28) KENDALL, M. G., STUART, A. (1979). *The advanced theory of statistics. Distributions theory*. London, Charles Griffin and Co., Ltd. 4$^{th}$ Edition.

29) KLEINROCK, L. (1985). *Queueing systems*. Vol. I and Vol. II. Wiley-New York.

30) RAMALHOTO, M. F., FERREIRA, M. A. M. (1996). Some further properties of the busy period of an M│G│∞ queue. *Central European Journal of Operations Research and Economics*, 4(4), 251-278.

31) ROSS, S. (1983). *Stochastic Processes*. Wiley, New York.

32) STADJE, W. (1985). The busy period of the queueing system M│G│∞. *Journal of Applied Probability*, 22, 697-704.

33) SYSKI, R. (1960). *Introduction to congestion theory in telephone systems*. Oliver and Boyd, London.

34) SYSKI, R. (1986). *Introduction to congestion theory in telephone systems*. North Holland, Amsterdam

35) TAKÁCS, L. (1962). *An introduction to queueing theory*. Oxford University Press, New York.

36) WHITT, W. (1984). On approximations for queues, I: extremal distributions. *A T &T Bell Laboratories Technical Journal*, 63(1), 115-138.